\documentclass{aomart}

\usepackage[english]{babel}
\usepackage{mathtools}

\def\Z{\mathbb{Z}}
\def\Q{\mathbb{Q}}
\def\R{\mathbb{R}}

\def\top{\mathsf{Top}}
\def\nbd{\mathsf{Nbd}}
\def\prtop{\mathsf{PreTop}}

\newcommand{\id}{\mathrm{id}}

\newtheorem*{prop}{Proposition}

\title{$\nbd$ is not cartesian closed}
\author{Giacomo Dossena}
\email{giacomo.dossena@gmail.com}
\copyrightnote{}
\thanks{April 2015}

\begin{document}

\begin{abstract}
In this short note it is shown that the category of neighborhood spaces is not cartesian closed.
\end{abstract}

\maketitle
\thispagestyle{plain}

A neighborhood space, or nbd space for short, is a set $X$ together with an assignment to each point $x\in X$ of a $p$-stack $\nu(x)$ on $X$ such that $\nu(x)\subset \dot{x}$. A $p$-stack on $X$ is just an upper collection of subsets of $X$ satisfying the pairwise intersection property (PIP): $A,B\in\nu(x)\implies A\cap B\neq\emptyset$. The formula $\nu(x)\subset \dot{x}$ means that each member of $\nu(x)$ contains the point $x$. Notice that if we insist that $\nu(x)$ be closed under finite intersections of its members then $\nu(x)$ becomes a filter (and the corresponding category becomes the category $\prtop$ of pretopological spaces).

A set map $f\colon X\to X'$ between ndb spaces $(X,\nu)$ and $(X',\nu')$ is said to be continuous at $x$ if $f(\nu(x))\supset \nu' (f(x))$, where $f(\nu(x))$ is the $p$-stack generated by the collection $\{f(A)\mid A\in\nu(x)\}$. In other words, the $p$-stack induced at $f(x)$ by transporting the $p$-stack at $x$ through $f$ refines the $p$-stack $\nu'(f(x))$. A set map $f\colon X\to X'$ is said to be continuous if it is continuous at each $x\in X$.

The collection of all ndb spaces and continuous maps between them forms a category, indicated by $\nbd$. The category $\nbd$ is topological with respect to the forgetful functor that assigns to each nbd space its underlying set. This means that each structured source has a unique initial lift (implying dually that each structured sink has a unique final lift). Moreover, the nbd spaces over a fixed set $X$ form a complete lattice with the partial order defined by $\nu\leq \nu'$ iff $\id_X\colon (X,\nu')\to (X,\nu)$ is continuous.

It is the purpose of this note to show that $\nbd$ is not cartesian closed. We will do that by constructing two quotient maps whose product is not quotient (in a cartesian closed topological category the product of two quotient maps is again quotient).

The maps are simply those of Example 4 in \cite{bhl1991} (used to show that $\top$ and $\prtop$ are not cartesian closed). Explicitly, let $\phi\colon \R\to \R/\Z$ be the map collapsing the integers to a point (so here $\R/\Z$ is the quotient space obtained by identifying all integers). This is a quotient map in $\top$ (by construction), but also in $\prtop$ (by explicit verification) and in $\nbd$ (either by explicit verification or by noticing that $\prtop$ embeds coreflectively in $\nbd$, so any quotient map in $\prtop$ is also quotient in $\nbd$). Let $\id_{\Q}\colon \Q\to \Q$ be the identity map on the rationals, obviously a quotient map in $\nbd$ (as well as in $\top$ and $\prtop$). Then form the product map $\Phi\coloneqq\phi\times\id_{\Q}\colon\R\times\Q\to\R/\Z\times\Q$. We claim that $\Phi$ is not quotient in $\nbd$ (actually \cite{bhl1991} shows that $\Phi$ is quotient neither in $\top$ nor in $\prtop$). To this end, we first recall the characterisation of quotient maps in $\nbd$ as given in \cite[Prop. 4.3]{km2002}.

\begin{prop}
Given a surjective set function $f\colon X\to Y$ and a nbd structure $\nu_X$ on $X$, the induced quotient nbd structure on $Y$ is given by $\nu_Y$ where, for each $y\in Y$, $\nu_Y(y)\coloneqq\{A\subset Y\mid \forall x\in f^{-1}(y), f^{-1}(A)\in\nu_X(x)\}$.
\end{prop}

Let us fix any $q\in\Q$ and define the following subsets of $\R\times\Q$:

\begin{equation*}
\begin{aligned}
A_q & \coloneqq \bigcup_{z\in\Z} \left(z-\tfrac{1}{2},z+\tfrac{1}{2}\right)\times \left(q-\tfrac{1}{1+|z|},q+\tfrac{1}{1+|z|}\right)\\
B_q & \coloneqq \Z\times(q-1,q+1)
\end{aligned}
\end{equation*}

Then $\Phi(A_q)$ does not belong to the neighborhood filter at $([0],q)$ (we indicate with $[0]$ the equivalence class of 0 in $\R/\Z$) because it does not contain any product of open sets in the component spaces $\R/\Z$ and $\Q$ (intuitively, this is because $A_q$ is made of boxes that shrink to smaller and smaller size in the $\Q$-component as $z$ moves away from 0). On the other hand, $\Phi(A_q)$ belongs to the $p$-stack at $([0],q)$ of the quotient nbd structure on $\R/\Z\times\Q$ induced by $\Phi$. To see this, notice that $\Phi^{-1}\Phi(A_q)=A_q\cup B_q$ belongs to the neighborhood filter at $(z,q)$ for each $z\in\Z$.

Notice that the argument carries over to the case where we substitute $\Q$ with $\R$. However, such a substitution destroys the validity of the example in $\top$ and $\prtop$ (so if we want a single example to prove non-cartesian-closedness of these three categories at once, it is better to stick to $\Q$).

\end{document}